\documentclass[12pt]{article}
\usepackage{amssymb, amsmath}

\newtheorem{theorem}{Theorem}[section]

\newcommand{\vare}{\varepsilon}
\newcommand{\n}{\nonumber}

\newcommand{\si}{\sigma_R }

\newcommand{\s}{\sigma}

\renewcommand{\o}{\omega}

\newcommand{\bb}{\begin{equation}}
\newcommand{\ee}{\end{equation}}
\newcommand{\bq}{\begin{eqnarray}}
\newcommand{\eq}{\end{eqnarray}}
\newcommand{\bqn}{\begin{eqnarray*}}
\newcommand{\eqn}{\end{eqnarray*}}

\begin{document}
\title{ On the uniqueness of solution to the steady Euler equations with perturbations}
\author{Dongho Chae\\
Department of Mathematics\\
            Chung-Ang University\\
          Seoul 156-756, Korea\\
              {\it e-mail : dchae@cau.ac.kr}}
 \date{}
\maketitle
\begin{abstract}
In this paper we study the uniqueness property of solutions to the
steady incompressible Euler equations with perturbations in $\Bbb
R^N$. Our perturbations include as special cases the Euler equations
with a `single signed' nonlinear term, the self-similar Euler
equations, and the steady Navier-Stokes equations. For these
equations show that suitable decay assumptions at infinity on the
solution or its derivatives, imposed by the $L^q$ conditions imply
that the only possible solution is zero. 
  \noindent{\bf
keywords:}   Euler equations with perturbation, steady solutions,
vanishing property
\end{abstract}
\section{ Main theorems}
 \setcounter{equation}{0}
We are concerned    on the  steady equations on $\Bbb R^N$ with
perturbation.
  \bb\label{main}
\left\{ \aligned & ( v\cdot \nabla )  v = -\nabla p +\Phi(v),\\
&\mathrm{div} \, v=0,
 \endaligned
 \right.
 \ee
 where $v=v(x)=(v_1(x),\cdots, v_N(x))$ is the velocity, and  $p=p(x)$ is the pressure.
  The function  $\Phi:\Bbb R^N\to \Bbb R^N$ defining the perturbation term satisfies suitable
 conditions depending on the cases we consider below.
  We study the vanishing property of the solutions to
 (\ref{main}).
 In this paper we consider the three cases of $\Phi(v)$. One
 is case where $\Phi(v)$ represents  a single signed nonlinear function(see below for more precise definition), and the
 other one is the case where the system (\ref{main}) corresponds to a
 generalization of the self-similar Euler equations, and finally the
 case where $\Phi(v)=\Delta v$, which corresponds to the steady
 Navier-Stokes equations.
\subsection{ The case where $\Phi(v)\cdot v$ is single signed}
Let us fix $N\geq 2$.  Here we assume that the continuous function
$\Phi(\cdot):\Bbb R^N\to \Bbb R^N$ satisfies the condition of single
signedness:
 \bb\label{da1}
\mbox{$\forall v\in \Bbb R^N$ either $ \Phi(v)\cdot v\geq 0$ or $
\Phi(v)\cdot v\leq 0$},
 \ee
   and
 \bb\label{da2}
 \mbox{ $\Phi(v)\cdot v=0$ if only if $v=0$}.
 \ee
   For such
given $\Phi$ we consider the system (\ref{main}).
 Note that when $\Phi(v)=-v$ the system (\ref{main})-(\ref{da2}) becomes the usual steady
 Euler equations with a damping term. More generally $\Phi(v)= G(|v|)v$ with a scalar function $G(x)>0$ for
$x>0$ satisfies (\ref{da1})-(\ref{da2}). We observe that the system
(\ref{main})-(\ref{da2}) has a trivial solution $v=0$. We will prove
that the uniqueness of solution to the system
(\ref{main})-(\ref{da2}) under quite mild decay conditions on the
solutions. More specifically we will prove the following.
\begin{theorem}
 Let $(v, p)$ be a $C^1(\Bbb R^N)$ solution of (\ref{main})-(\ref{da2}). Suppose  there exists $q\in
[\frac{3N}{N-1} , \infty)$ such that
 \bb\label{13a}
   v\in L^q (\Bbb R^N) \quad \mbox{and}\quad p\in L^{\frac{q}{2}}
   (\Bbb R^N).
  \ee
  Then, $v=0$.
  \end{theorem}
{\em Remark 1.1 } If $\Phi(v)$ satisfies an extra condition
div$\Phi(v(x))=0$, then we do not need to assume $p\in
L^{\frac{q}{2}}
   (\Bbb R^N)$ in (\ref{13a}). Since in that case we have the well-known
   velocity-pressure relation as in the incompressible
   Euler or the Navier-Stokes equations,
   $$p(x)=\sum_{j,k=1}^N R_j R_k (v_jv_k)(x) $$
    with the Riesz transform $R_j$, $j=1,\cdots,N,$ in $\Bbb R^N$
    (\cite{ste}),
and the $L^{\frac{q}{2}}$ estimate of the pressure follows from the
$L^q$ estimate for the velocity by the Calderon-Zygmund inequality,
 \bb\label{cz}\|p\|_{L^{\frac{q}{2}}}\leq C\sum_{j,k=1}^N  \|R_jR_k v_jv_k \|_{L^{\frac{q}{2}}}\leq C
 \|v\|_{L^{q}}^2 \quad 2<q<\infty.
 \ee

\subsection{ The case $\Phi(v)=av +b(x\cdot \nabla) v, \, ab\neq0$}
In this subsection we fix $N=3$. Let $a,b$ are given constants such
that $ab\neq 0$. We study here the system in $\Bbb R^3$.
 \bb\label{pe}
\left\{ \aligned & ( v\cdot \nabla )  v = -\nabla p  +a v+b(x\cdot \nabla) v,\\
&\mathrm{div} \, v=0.
 \endaligned
 \right.
 \ee
In the special case of $a=-\frac{\alpha}{\alpha+1}$,
$b=-\frac{1}{\alpha+1}$ the system (\ref{pe}) reduces to the
self-similar Euler equations.
 \bb\label{ss}
\left\{ \aligned & \frac{\alpha}{\alpha+1} v+\frac{1}{\alpha+1}(x\cdot \nabla) v +( v\cdot \nabla )  v = -\nabla p ,\\
&\mathrm{div} \, v=0.
 \endaligned
 \right.
 \ee
The system (\ref{ss}) is obtained from the time dependent Euler
equations,
 $$
(E)\left\{\aligned &\partial_t u +( u\cdot \nabla )  u = -\nabla \mathbf{p} ,\\
&\mathrm{div} \, u=0,
 \endaligned\right.
 $$
by the self-similar ansatz,
 \bqn
u(x,t)&=&\frac{1}{(t_*-t)^{\frac{\alpha}{\alpha+1}}}
v\left(\frac{x-x_*}{(t_*-t)^{\frac{1}{\alpha+1}}} \right),\\
\mathbf{p}(x,t)&=&\frac{1}{(t_*-t)^{\frac{2\alpha}{\alpha+1}}}
p\left(\frac{x-x_*}{(t_*-t)^{\frac{1}{\alpha+1}}} \right),
 \eqn
 where $(x_*,t_*)$ is the hypothetical self-similar
 blow-up space-time point.
 The question of self-similar blow-up for the Navier-Stokes
 equations is asked in (\cite{ler}), and is answered negatively
 in \cite{nec} for $v\in L^3(\Bbb R^3)$, and is extended in \cite{tsa} for $v\in L^q(\Bbb R^3), q\geq 3$. Similar problem for the Euler equations is
 studied in \cite{cha0, cha00}.
 For $\alpha <\infty$ with $\alpha \neq -1$
it is proved in \cite{cha0} that if  a solution to (\ref{ss}), $v\in
C^1(\Bbb R^3)$, decaying to zero at infinity, satisfies
$\o=\mathrm{curl}\, v\in \cap_{0<q<q_0} L^q(\Bbb R^3)$ for some
$q_0>0$, then $v=0$. In the extreme case
 $\alpha=\infty$, we have (\ref{ss}) becomes the system (\ref{main}) with
 $\Phi(v)=-v$.
\begin{theorem} Let $v$ be a classical solution to (\ref{pe}). Suppose
there exists $q_0>0$ such that
 \bb\label{condi13} \|\nabla v\|_{L^\infty}<\infty,
\quad\mbox{and}\quad \o\in \bigcap_{0<q< q_0} L^q(\Bbb R^3).
  \ee
 Then, $v=\nabla h$ for a  harmonic scalar function $h$
on $\Bbb R^3$. Thus, if we impose further the condition $
\lim_{|x|\to \infty} |v(x)|=0$, then $v=0$.
\end{theorem}
{\em Remark 1.2 } In \cite{cha0} we used the time dependent Euler
equations directly to prove Theorem 1.1, and needed the decay
condition for the velocity, since we used the notion of
back-to-label map, whose existence is guaranteed for the decaying
velocity(\cite{con0}). In the proof of the above theorem below,
however, we work with the stationary system (\ref{ss}), and do not
use the back-to-label map, and therefore the
decay condition for the velocity field is not necessary.\\
\ \\
\noindent{\em Remark 1.3 } As far as the regularity assumption for
the solution $v$,  what we need in the proof  is actually the
differentiability almost
everywhere, which is guaranteed by the first condition of (\ref{condi13}).\\

\subsection{ The case $\Phi(v)=a\Delta v, \, a\neq0$}
In this subsection we also fix $N=3$.  Here we study (\ref{main})
with $\Phi(v)=a\Delta v$. In this case without loss of generality we
may set $a=1$. In this case the system (\ref{main}) reduces to the
steady Navier-Stokes equations in $\Bbb R^3$.
 $$
(NS)\left\{ \aligned & ( v\cdot \nabla )  v = -\nabla p +\Delta v,\\
&\mathrm{div} \, v=0,
 \endaligned
 \right.
 $$
We consider here the generalized solutions of the system (NS),
satisfying
 \bb\label{diri}
\int_{\Bbb R^3} |\nabla v|^2dx <\infty,
 \ee
 and
    \bb\label{12c}
  \lim_{|x|\to \infty} v(x)= 0.
  \ee
 It is well-known that a generalized solution to (NS) belonging to $ W^{1,2}_{loc}(\Bbb R^3)$ implies that
 $v$ is smooth(see e.g.\cite{gal}). Therefore without loss of
 generality we can assume that our solutions to (NS) satisfying (\ref{diri}) are smooth.
 The uniqueness question, or equivalently the question of  Liouville
 property of solution for the system (NS) under the assumptions
 (\ref{diri}) and (\ref{12c}) is a long standing open problem.
  On the other hand, it is well-known that the uniqueness of solution holds in the
   class $L^{\frac92} (\Bbb R^3)\cap \dot{H}^1(\Bbb R^3)$, namely a
  smooth solution to (NS) satisfying
$v\in L^{\frac92} (\Bbb R^3)$
   and (\ref{diri}) is $v=0$(see Theorem 9.7 of \cite{gal}). We assume here slightly stronger condition than
  (\ref{diri}), but having the same scaling property, to deduce the uniqueness result. More precisely, we
  have the following theorem.
\begin{theorem}
Let $v$ be a smooth solution of (NS) satisfying (\ref{12c}) and
 \bb\label{sdiri}
 \int_{\Bbb R^3} |\Delta v|^{\frac65}\, dx<\infty.
 \ee
   Then, $v=0$ on
 $\Bbb R^3$.
  \end{theorem}
  {\em Remark 1.3 } Under the assumption (\ref{12c}) we have the
  inequalities with the norms of the {\em same scaling properties,}
$$
 \|\nabla v\|_{L^2} \leq C  \|D^2 v\|_{L^{\frac65}}
 \leq C\|\Delta v\|_{L^\frac65}<\infty
$$
due to the Sobolev and the Calderon-Zygmund inequalities. Thus,
(\ref{sdiri}) implies (\ref{diri}).  There is no, however, mutual
implication relation between Theorem 1.3 and
the above mentioned $L^{\frac92}$ result.\\
\ \\
{\em This paper is a modified and extended version of author's
preprint \cite{cha000}.}
\section{Proof of the Main Theorems }
\setcounter{equation}{0}

\noindent{\bf Proof of Theorem 1.1 }  We denote
 $$[f]_+=\max\{0, f\}, \quad  [f]_-=\max\{0, -f\},$$
and
$$
D_\pm:=\left\{ x\in \Bbb R^N\, \Big|\,\left[p(x)+\frac12
|v(x)|^2\right]_\pm >0\right\}
$$
respectively. We introduce the radial cut-off function $\sigma\in
C_0 ^\infty(\Bbb R^N)$ such that
 \bb\label{16}
   \sigma(|x|)=\left\{ \aligned
                  &1 \quad\mbox{if $|x|<1$},\\
                     &0 \quad\mbox{if $|x|>2$},
                      \endaligned \right.
 \ee
and $0\leq \sigma  (x)\leq 1$ for $1<|x|<2$.  Then, for each $R
>0$, we define
 $$
\s \left(\frac{|x|}{R}\right):=\s_R (|x|)\in C_0 ^\infty (\Bbb R^N).
$$
  We multiply first equations of (\ref{main}) by $v$ to obtain
   \bb\label{main1}
   v\cdot \Phi(v) =v\cdot \nabla \left(p+\frac12 |v|^2\right).
   \ee
   Next, we multiply (\ref{main1}) by
  $ \left[p+\frac12 |v|^2\right]_+
^{\frac{qN-q-3N}{2N}}\s_R \, \mathrm{sign}\{v\cdot \Phi(v)\} $ and
integrate over $\Bbb R^N$, then we have
  \bq\label{euler1}
 \lefteqn{\int_{\Bbb R^N} \left[p+\frac12
|v|^2\right]_+^{\frac{qN-q-3N}{2N}}\left|v\cdot \Phi(v)\right|
\s_R\,
dx}\n \\
&&=\mathrm{sign}\{ v\cdot \Phi(v)\}\int_{\Bbb R^N} \left[p+\frac12
|v|^2\right]_+^{\frac{qN-q-3N}{2N}}\s_R v \cdot\nabla
\left(p+\frac12 |v|^2\right) \,dx\n \\
&&:=I
   \eq
   We estimate $I$ as follows.
 \bqn
 |I|&=&\left|\int_{\Bbb R^N} \left[p+\frac12 |v|^2\right]_+^{\frac{qN-q-3N}{2N}} \s_R
v\cdot\nabla \left(p+\frac12 |v|^2 \right)\, dx\right| \n \\
&=&\left|\int_{D_+} \left[p+\frac12
|v|^2\right]_+^{\frac{qN-q-3N}{2N}} \s_R
v\cdot\nabla \left[p+\frac12 |v|^2 \right]_+\, dx\right|\n \\
 &=& \frac{2N}{qN-q-N} \left|\int_{D_+}\si v\cdot\nabla \left[p+\frac12
|v|^2 \right]_+^{\frac{qN-q-N}{2N}} \, dx \right|\n \\
&=&\frac{2N}{qN-q-N} \left| \int_{D_+}\left[p+\frac12 |v|^2
\right]_+^{\frac{qN-q-N}{2N}}  v\cdot\nabla \si \, dx\right| \eqn
\bq \label{n116}
 &\leq&
\frac{C\|\nabla\s\|_{L^\infty} }{R} \left(\int_{\Bbb R^N}
(|p|+|v|^2)^{\frac{q}{2}} \, dx\right)^{\frac{qN-q-N}{qN}}
\|v\|_{L^q(R\leq |x|\leq 2R)} \times\n \\
&&\hspace{.5in} \times\left(\int_{\{ R\leq |x|\leq 2R\}} \,
dx \right)^{\frac1N}\n \\
  &\leq& C\|\nabla\s\|_{L^\infty}
\left(\|p\|_{L^{\frac{q}{2}}} +\|v\|_{L^q}^2
\right)^{\frac{qN-q-N}{qN}}\|v\|_{L^q(R\leq |x|\leq 2R)}\to 0
 \eq
as $R\to \infty$.
 Passing $R\to \infty$ in (\ref{euler1}), we obtain
 \bb\label{116a}
 \int_{\Bbb R^N} \left[p+\frac12
|v|^2\right]_+^{\frac{qN-q-3N}{2N}} \left|v\cdot \Phi(v)\right| \,
dx =0. \ee
  Similarly, multiplying (\ref{main1}) by $ \left[p+\frac12
  |v|^2\right]_-
^{\frac{qN-q-3N}{2N}} \s_R $, and integrate over $\Bbb R^N$,  we
deduce by similar estimates to the above,
 \bq\label{116aa}
 \lefteqn{\int_{\Bbb R^N}
 \left[p+\frac12 |v|^2\right]_-^{\frac{qN-q-3N}{2N}} \left|v\cdot \Phi(v)\right| \s_R\, dx}\hspace{.3in}\n \\
 &&=-\int_{\Bbb R^N} \left[p+\frac12 |v|^2\right]_-^{\frac{qN-q-3N}{2N}} \s_R
v\cdot\nabla \left(p+\frac12 |v|^2 \right)\, dx \n \\
&&=\int_{\Bbb R^N} \left[p+\frac12
|v|^2\right]_-^{\frac{qN-q-3N}{2N}} \s_R
v\cdot\nabla \left[p+\frac12 |v|^2 \right]_-\, dx\n \\
&&\leq C\|\nabla\s\|_{L^\infty}\left(\|p\|_{L^{\frac{q}{2}}}
+\|v\|_{L^q}^2 \right)^{\frac{qN-q-N}{qN}}\|v\|_{L^q(R\leq |x|\leq
2R)}\to 0\n \\
 \eq
as $R\to \infty$. Hence,
 \bb\label{116ab}\int_{\Bbb R^N}
\left[p+\frac12 |v|^2\right]_-^{\frac{qN-q-3N}{2N}}\left|v\cdot
\Phi(v)\right| \, dx =0. \ee Let us define
$$ \mathcal{S}=\{ x\in \Bbb R^N\, |\, v(x)\neq0\}. $$
Suppose $\mathcal{S}\neq \emptyset$. Then, (\ref{116aa}) and
(\ref{116ab}) together with (\ref{da1})-(\ref{da2}) imply
$$
\left[p(x)+\frac12 |v(x)|^2\right]_+=\left[p(x)+\frac12
|v(x)|^2\right]_-=0\quad \forall x\in \mathcal{S}.
$$
Namely,
$$
p(x)+\frac12 |v(x)|^2=0\quad \forall x\in \mathcal{S}.
$$
From (\ref{main1}) this implies
 \bb\label{119}
  \Phi(v(x))\cdot v(x)=0\qquad \forall x\in \mathcal{S}.
  \ee
 Considering the conditions on $\Phi$ in (\ref{da1})-(\ref{da2}), we have a contradiction, and therefore we need $\mathcal{S}=\emptyset$, namely
 $v=0$ on $\Bbb R^N$.
   $\square$\\
\ \\
 \noindent{\bf Proof of Theorem 1.2 }
 We first observe that from the calculus identity
 $$ v(x)=v(0)+\int_0 ^1\partial_s v(sx) ds=v(0)+\int_0 ^1 x\cdot
 \nabla v(sx) ds,
 $$
 we have
 $|v(x)|\leq |v(0)|+ |x|\|\nabla v\|_{L^\infty}\leq C(1+|x|)\|\nabla
 v\|_{L^\infty},$
 and
  \bb\label{cal}
  \sup_{x\in \Bbb R^3} \frac{|v(x)|}{1+|x|} \leq C\|\nabla
  v\|_{L^\infty}.
  \ee
We consider the vorticity equation of (\ref{pe}),
 \bb\label{vor}
 -(a+b)\o- b( x\cdot \nabla ) \o+(v\cdot \nabla
 )\o=(\o\cdot \nabla )v .
 \ee
 We take $L^2(\Bbb R^3)$ inner product (\ref{vor}) by $|\o|^{q-2}\o
 \si$, then after integration by part, we have
 \bq\label{sel}
 \lefteqn{-(a+b)\|\o\si\|_{L^q}^q+\frac{3b}{q} \|\o\si\|_{L^q}^q-\int_{\Bbb R^3}
\xi\cdot \nabla v\cdot \xi |\o|^q \si\, dx}\n \\
&&= b \int_{\Bbb R^3} |\o|^q (x\cdot\nabla )\si \,
dx +\int_{\Bbb R^3} |\o|^q (v\cdot\nabla )\si \, dx\n \\
&&:=I+J.
 \eq
We estimate $I$ and $J$ easily as follows.
  $$
  |I|\leq  \frac{b}{ R}\int_{\{R\leq |x|\leq 2R\}} |\o|^q
  |x||\nabla \s|\, dx
  \leq b\|\nabla \s\|_{L^\infty}\|\o\|^q_{L^p(R\leq |x|\leq
  2R)}\to 0
  $$
   as $R\to \infty$.
   \bqn
    |J|&\leq&  \frac{1}{ R}\int_{\{R\leq |x|\leq 2R\}} |\o|^q
  |v||\nabla \s|\, dx
  \leq\frac{1+2R}{ R}\int_{\{R\leq |x|\leq 2R\}}
  \frac{|v(x)|}{1+|x|}|\o|^q|\nabla \s|\, dx\n \\
  &\leq &\frac{1+2R}{ R}\|\nabla \s\|_{L^\infty}\|\nabla v\|_{L^\infty}\|\o\|^q_{L^p(R\leq |x|\leq
  2R)}\to 0
  \eqn
  as $R\to \infty$, where we used (\ref{cal}). Therefore, passing
  $R\to \infty$ in (\ref{sel}), we obtain
  \bb\label{sel2} -\|\nabla v\|_{L^\infty}\|\o\|_{L^q}^q
\leq \left(a+b-\frac{3b}{q}\right)\|\o\|_{L^q}^q\leq \|\nabla
v\|_{L^\infty}\|\o\|_{L^q}^q.
  \ee
 Suppose  $\o\neq 0$, then we will derive a contradiction. If $\|\o\|_{L^q} \neq 0$, we can divide (\ref{sel2}) by
 $\|\o\|_{L^q}^q$ to have
\bb\label{sel3}
 -\|\nabla v\|_{L^\infty} \leq
\left(a+b-\frac{3b}{q}\right)\leq \|\nabla v\|_{L^\infty},
  \ee
 which holds for all $q\in (0, q_0)$. Since $b\neq 0$, passing $q\downarrow 0$ in (\ref{sel3}), we
 obtain desired contradiction. Therefore $\o=\mathrm{curl}\, v=0$. This,
 together with $\mathrm{div}\,v=0$, provides us with the fact that $v=\nabla
 h$ for a scalar harmonic function $h$ on $\Bbb R^3$. $\square$\\
 \ \\
\noindent{\bf Proof of Theorem 1.3}
 Under the
assumption (\ref{sdiri}) and Remark 1.1, Theorem IX.6.1 of
\cite{gal} implies
 that
  \bb\label{decay}
 \lim_{|x|\to \infty}|p(x)-p_1 |=0.
  \ee
 for  a constant $p_1$.
 Therefore, if we set
 $$
  \label{ber}Q(x):=\frac12 |v(x)|^2 +p(x)-p_1,
 $$
 then
  \bb\label{13}
\lim_{|x|\to \infty} |Q(x)|=0.
  \ee
As before we  denote $[f]_+=\max\{0, f\}, \quad  [f]_-=\max\{0,
-f\}.$ Given $\vare > 0$, we define
  \bqn
    D_+^\vare&=&\left\{ x\in \Bbb
R^3\, \Big|\,\left[Q(x)-\vare\right]_+>0\right\},\n
\\
D_-^\vare&=&\left\{ x\in \Bbb R^3\,
\Big|\,\left[Q(x)+\vare\right]_->0\right\}.
  \eqn
 respectively. Note
that (\ref{13}) implies that $D_\pm^\vare$ are bounded sets in $\Bbb
R^3$. Moreover,
  \bb\label{obs} Q\mp\vare
=0\quad\mbox{on}\quad \partial D_\pm^\vare
 \ee
 respectively.
Also, thanks to the Sard theorem combined with the implicit function
theorem $\partial D_\pm^\vare$'s are smooth level surfaces in $\Bbb
R^3$ except the values of $\vare>0$, having the zero Lebesgue
measure, which corresponds to the critical values of $z=Q(x)$. It is
understood that our values of  $\vare$ below avoids these
exceptional ones. We write the system (NS) in the form,
 \bb\label{17}
 -v\times\mathrm{ curl} \, v =-\nabla Q +\Delta v.
 \ee
 Let us
multiply (\ref{17}) by $ v \left[Q-\vare\right]_+$, and integrate it
over $\Bbb R^3$. Then, since
  $
v\times \mathrm{curl}\,v  \cdot v =0, $
  we have
 \bq\label{18}
0&=& -\int_{\Bbb R^3} \left[Q-\vare\right]_+ v \cdot\nabla
\left(Q-\vare\right) \,dx+\int_{\Bbb R^3}\left[Q-\vare\right]_+ v\cdot \Delta v\, dx\n\\
&:=&I_1 +I_2 .
 \eq
Integrating by parts, using (\ref{obs}), we obtain
 \bqn
 I_1=-\int_{D_+^\vare} \left(Q-\vare\right)
v\cdot\nabla \left(Q -\vare\right)\, dx= -\frac{1}{2}
\int_{D_+^\vare}
  v\cdot\nabla \left(Q-\vare \right)^{2} \, dx.
=0 \eqn
 Using
 \bb\label{f2aa}
 v\cdot \Delta v=\Delta (\frac12|v|^2)-|\nabla v|^2,
 \ee
 and the well-known formula for the Navier-Stokes equations,
 \bb\label{f2a}
 \Delta p=|\o|^2-|\nabla v |^2,
 \ee
 we have
 \bq\label{110}
 I_2&=&-\int_{\Bbb R^3} |\nabla v|^2\left[ Q-\vare\right]_+ \, dx +\int_{\Bbb R^3}\Delta \left(\frac12 |v|^2\right)
\left[Q-\vare\right]_+ \, dx \n \\
  &=&-\int_{\Bbb R^3}|\o|^2\left[Q-\vare\right]_+\, dx +\int_{\Bbb R^3} \Delta \left(Q-\vare\right)
 \left[ Q-\vare\right]_+
 \, dx\n \\
&:=&J_1+J_2.
 \eq
Integrating by parts, we transform $J_2$ into
 \bq\label{110a}
 J_2=\int_{D_+^\vare} \Delta\left(Q-\vare\right)
 \left( Q-\vare\right)
 \,dx=-
 \int_{D_+^\vare}\left|\nabla\left(Q-\vare\right)\right|^2
 \,dx.
 \eq
Thus,  the derivations (\ref{18})-(\ref{110a})
 lead us to
 \bb\label{dom}
0=\int_{D_+^\vare} |\o|^2\left| Q-\vare\right|\,
dx+\int_{D_+^\vare}\left|\nabla\left(Q-\vare\right)\right|^2
 \,dx
\ee for all $\vare>0$. The vanishing of the second term of
(\ref{dom}) implies
$$
\left[Q-\vare\right]_+=C_0\quad \mbox{on}\quad D_+^\vare
$$ for a constant $C_0$. From the fact (\ref{obs})
 we have $C_0=0$, and $[Q-\vare]_+=0$ on $\Bbb R^3$, which holds for all  $\vare >0.$
Hence,
 \bb\label{plus}
\left[Q\right]_+=0 \quad\mbox{on $\Bbb R^3$}.
  \ee
This shows that $Q\leq 0$ on $\Bbb R^3$. Suppose $Q=0$ on $\Bbb
R^3$. Then, from (\ref{17}), we have $v\cdot\Delta v=0$ on $\Bbb
R^3$. Hence,
$$ \Delta p=-\frac12\Delta |v|^2=-v\cdot\Delta v-|\nabla v|^2=-|\nabla
v|^2.
$$
Comparing this with (\ref{f2a}), we have $\o=0$. Combining this with
div $v=0$, we find that $v$ is a harmonic function in $\Bbb R^3$.
Thus, by (\ref{12c}) and the Liouville theorem for the harmonic
function, $v=0$. Hence, without loss of generality, we may assume
$$0>\inf_{x\in \Bbb R^3} Q(x):=-\delta_0.
$$
Given $\delta \in (0, \delta_0)$, we multiply (\ref{17}) by $ v
\left[Q+\vare\right]_- ^{\delta} $, and integrate it over $\Bbb
R^3$. Then, similarly to the above
  we have
 \bq\label{m18}
0&=& -\int_{\Bbb R^3} \left[Q+\vare\right]_-^{\delta} v \cdot\nabla
\left(Q+\vare\right) \,dx+\int_{\Bbb R^3} \left[Q+\vare\right]_-^{\delta} v\cdot \Delta v\, dx\n\\
&:=&I_1' +I_2' .
 \eq
Observing $Q(x)+\vare=- \left[Q(x)+\vare\right]_-$ for all $x\in
D_-^\vare$, integrating by part, we obtain
 \bqn
 I_1'&=&\int_{D_-^\vare}\left[Q+\vare\right]_-^{\delta}
v\cdot\nabla \left[Q+\vare\right]_-\, dx\n \\
 &=&\frac{1}{\delta+1} \int_{D_-^\vare} v\cdot\nabla \left[Q+\vare \right]_-^{\delta+1} \, dx =0.
 \eqn
 Thus, using (\ref{f2aa}), we have
 \bq\label{m110}
 0=-\int_{D_-^\vare} |\nabla v|^2\left[Q+\vare\right]_-^{\delta} \,
 dx+\frac12 \int_{D_-^\vare}\left[ Q+\vare\right]_- ^{\delta}\Delta |v|^2
\, dx \eq
 Now,  we have the point-wise convergence
 $$\left[Q(x)+\vare\right]_-^\delta \to 1 \quad \forall x\in D_-^\vare.
 $$
 as $\delta\downarrow 0$. Hence, passing $\delta\downarrow 0$ in
 (\ref{m110}),
 by the dominated convergence
 theorem, we obtain
 \bb\label{m110a}
 \int_{D_-^\vare} |\nabla v|^2\,dx= \frac12\int_{D_-^\vare}
  \Delta |v|^2\, dx,
 \ee
which holds for all $\vare>0$. For a sequence $\{\vare_n\}$ with
$\vare_n \downarrow 0$ as $n\to \infty$, we observe
 $$ D_-^{\vare_n }\uparrow \cup_{n=1}^\infty D_-^{\vare_n}=D_-:=
 \{ x\in \Bbb R^3\, |\, Q(x)<0.\}.
 $$
Since
 \bqn
\left|\int_{\Bbb R^3} v\cdot \Delta v\, dx\right|&\leq&
\|v\|_{L^6}\|\Delta v\|_{L^{\frac65}}\leq C \|\nabla
v\|_{L^2}\|\Delta v\|_{L^{\frac65}}\n \\
&\leq &C \|\Delta v\|_{L^{\frac65}}^2<\infty,
 \eqn
 we have
 \bb\label{elone}\Delta |v|^2=2 v\cdot \Delta v +2 |\nabla v|^2 \in L^1(\Bbb R^2).
 \ee

Thus,  we can apply the dominated convergence theorem in passing
$\vare \downarrow 0$ in (\ref{m110a}) to deduce
 \bb\label{mmdom1}
\int_{D_-} |\nabla v|^2\, dx=\frac12\int_{D_-}\Delta |v|^2\, dx.
 \ee
 Now, thanks to (\ref{plus}) the set
 $$ S=\{ x\in \Bbb R^3\, |\, Q(x)=0\} $$
 consists of critical(maximum) points of $Q$, and hence
 $ \nabla Q(x)=0$ for all $x\in S,$ and the system (\ref{17}) reduces to
  \bb\label{redu}
   -v\times \o=\Delta v \quad\mbox{on}\quad S.
  \ee
Multiplying (\ref{redu}) by $v$, we have that
$$ 0=v\cdot \Delta v=\frac12\Delta |v|^2-|\nabla v|^2\quad\mbox{on}\quad S.
$$
Therefore, one can extend the domain of integration in
(\ref{mmdom1}) from $D_-$ to $D_-\cup S=\Bbb R^3$, and therefore
 \bb\label{mmdom2}
\int_{\Bbb R^3} |\nabla v|^2\, dx=\frac12\int_{\Bbb R^3}\Delta
|v|^2\, dx.
 \ee
 We now claim the right hand side of (\ref{mmdom2}) vanishes.
 Since $\Delta |v|^2\in L^1(\Bbb R^3)$ from (\ref{elone}), applying the dominated
 convergence theorem, we have
 \bqn
 \left|\int_{\Bbb R^3} \Delta |v|^2\, dx\right|&=&\lim_{R\to \infty} \left|\int_{\Bbb R^3}
 \Delta |v|^2 \si \, dx\right|
 =\lim_{R\to \infty} \left|\int_{\Bbb R^3}
  |v|^2 \Delta\si \, dx\right|\n \\
  &\leq&\lim_{R\to \infty}\int_{\Bbb R^3}
  |v|^2 |\Delta\si| \, dx\n \\
 &\leq& \lim_{R\to \infty}\frac{\|D^2\s\|_{L^\infty}}{R^2}
 \|v\|_{L^6(R\leq |x|\leq
 2R)}^2
 \left(\int_{\{R\leq |x|\leq
 2R\}} dx\right)^{\frac23}\n \\
 &\leq &C\|D^2\s\|_{L^\infty}\lim_{R\to \infty}
 \|v\|_{L^6(R\leq |x|\leq
 2R)}^2=0
 \eqn
 as claimed.
 Thus (\ref{mmdom2}) implies that
 $$
\nabla v=0 \quad\mbox{on} \quad \Bbb R^3,
$$
and $v=$ constant. By (\ref{12c}) we have $v=0$.
$\square$\\
\ \\
 \noindent{\em Remark after the proof of Theorem 1.3: } The first part of the
above proof, showing $ [Q]_+=0 $ can be also done by applying the
maximum principle, which is from the following identity for $Q$,
$$ -\Delta Q+v\cdot \nabla Q =-|\o|^2 \leq 0
$$
I do not think, however, the maximum principle can also be applied
to the proof of the second part, showing $[Q]_-=0$, which is more
subtle than the first part. The above proof overall shows that the
argument of the proof I used for this second part can also be
adapted for the first part without using the
maximum principle, which exhibits consistency.\\

 $$ \mbox{\bf Acknowledgements} $$
This work was supported partially by the NRF grant. no. 2006-0093854
and also by Chung-Ang University Research Grants in 2012.

\end{document}